\documentclass[11pt]{article}

\usepackage{amsmath,amsthm,amstext,amscd,amssymb,euscript,url}
\usepackage{mathrsfs}
\usepackage{calrsfs}
\usepackage{epsfig}
\usepackage[inline]{enumitem}

\newcommand{\Z}{\mathbb Z}
\newcommand{\R}{\mathbb R}

\newcommand{\E}{\mathbb E}
\newcommand{\Zd}{\mathbb Z^d}

\renewcommand{\phi}{\varphi}
\newcommand{\epsi}{\ensuremath{\epsilon}}

\newcommand{\si}{\ensuremath{\sigma}}

\newcommand{\pee}{\ensuremath{\mathbb{P}}}

\newcommand{\loc}{\mathcal{L}}

\def\1{{\mathchoice {\rm 1\mskip-4mu l} {\rm 1\mskip-4mu l}
{\rm 1\mskip-4.5mu l} {\rm 1\mskip-5mu l}}}

\newtheorem{theorem}{{\small T}{\scriptsize HEOREM}}[section]
\newtheorem{corollary}{{\bf{\small C}{\scriptsize OROLLARY}}}[section]
\newtheorem{proposition}{{\bf{\small P}{\scriptsize ROPOSITION}}}[section]
\newtheorem{lemma}{{\bf{\small L}{\scriptsize EMMA}}}[section]
\newtheorem{remark}{{\bf{\small R}{\scriptsize EMARK}}}[section]
\newtheorem{definition}{{\bf{\small D}{\scriptsize EFINITION}}}[section]

\renewenvironment{proof}[1]
{\noindent{{\bf{\small{ P}{\scriptsize ROOF}}}.}\hspace{0.1cm} #1} {$\;\qed$\newline}

\newcommand{\beq}{\begin{eqnarray}}
\newcommand{\eeq}{\end{eqnarray}}

\newcommand{\ba}{\begin{align*}}
\newcommand{\ea}{\end{align*}}

\newcommand{\be}{\begin{equation}}
\newcommand{\ee}{\end{equation}}

\newcommand{\bl}{\begin{lemma}}
\newcommand{\el}{\end{lemma}}

\newcommand{\br}{\begin{remark}}
\newcommand{\er}{\end{remark}}

\newcommand{\bt}{\begin{theorem}}
\newcommand{\et}{\end{theorem}}

\newcommand{\bd}{\begin{definition}}
\newcommand{\ed}{\end{definition}}

\newcommand{\bp}{\begin{proposition}}
\newcommand{\ep}{\end{proposition}}

\newcommand{\bc}{\begin{corollary}}
\newcommand{\ec}{\end{corollary}}

\newcommand{\bpr}{\begin{proof}}
\newcommand{\epr}{\end{proof}}

\newcommand{\bi}{\begin{itemize}}
\newcommand{\ei}{\end{itemize}}

\newcommand{\ben}{\begin{enumerate}}
\newcommand{\een}{\end{enumerate}}

\newcommand{\caF}{{\mathcal F}}

\newcommand{\caM}{{\mathcal M}}
\newcommand{\caN}{{\mathcal N}}

\newcommand{\caP}{{\mathcal P}}

\newcommand{\caS}{{\mathcal S}}

\newcommand{\caV}{{\mathcal V}}

\newcommand{\caX}{{\mathcal X}}

\newcommand{\omu}{\overline{\mu}}
\newcommand{\homu}{\widehat{\mu}}

\begin{document}
\title{{\bf Central limit theorem and large deviations for run and tumble particles on a lattice}}
\author{Bart van Gisbergen and
Frank Redig\\
Technische Universiteit Delft\\ Van Mourik Broekmanweg 6,\\ 2628 XE, Delft, The Netherlands}
\maketitle
\begin{abstract}
We study  run and tumble  particles on the one-dimensional lattice $\Z$.
We explicitly compute the Fourier-Laplace transform of the position of the particle and as a consequence obtain explicit expressions for the diffusion constant and the large deviation free energy function.
We also do the same computations in a corresponding continuum model. In the latter, when adding an external field, we can
explicitly compute the large deviation free energy, and the deviation from the Einstein relation
due to activity. Finally, we generalize the model to the $d$-dimensional lattice $\Zd$, with an arbitrary finite set of velocities, and show that the large deviation free energy for the position of the particle can be computed via the largest eigenvalue of a matrix of
Schr\"{o}dinger operator form, for which we can derive an explicit variational formula via occupation time large deviations of the velocity flip process.
\end{abstract}

\section{Introduction}
Run and tumble particles are simple models of active matter, where particles move under influence of an internal degree of freedom, and external noise
\cite{maes, majumdar, dhar, dhar2}.
The activity of the particles is a source of non-equilibrium and it is of interest to understand macroscopic behaviour as a function of the  parameters modelling the activity of the particles. Contrary to an external driving field, the effects of activity typically appear as second order effects (for small activity).
In simple lattice (and continuum) models, we can quantitatively understand the influence of activity on transport coefficients, and large deviation functions, and compare with the equilibrium, or close to equilibrium setting.
In this paper we first consider simple models of such motion on a lattice. This has the advantage that we can do explicit computations of all the relevant macroscopic quantities: diffusion constant, and large deviation free energy function. Notice that in the continuum, large deviation results have been obtained before in \cite{bombo}.
Finally we analyze a general model on $\Zd$ with a finite set of velocities, and obtain a variational formula for the large deviation free energy
of the particle's position.
The rest of our paper is organized as follows.
In section 2 we consider the basic lattice model, compute its characteristic function, and prove the central limit theorem for the particle's position. We also consider a continuum limit, where we can do the same computations, recovering the results in \cite{bombo}.
In section 3 we study the large deviations for the position of the particle by explicitly computing the large deviation free energy function. We show that the continuum limit of the
large deviation free energy converges to the large deviation free energy of the continuum model, previously computed in \cite{bombo}.
In the continuum model with drift, we can explicitly compute the deviation from the Einstein relation due to activity.
In section 4 we consider various generalizations of the model and prove a connection between the large deviation free energy of the particle's position and the occupation time
large deviation function for the autonomous velocity flip process.

\section{The model and its scaling behavior}
\subsection{The model}
The active particle has a position $x\in\Z$ and a velocity $v\in \{-1,1\}$.
The process $\{(X_t,v_t): t\geq 0\}$ is described via the generator
\beq\label{gena}
Lf(x,v) &=&  \lambda (f(x+v,v)- f(x,v)) \nonumber\\
&+& \kappa (f(x+1, v)+ f(x-1,v) -2f(x,v))
\nonumber\\
&+ & \gamma(f(x,-v)-f(x,v))
\eeq
This is interpreted as follows: with rate $\lambda$ the process makes a jump in the direction of the velocity,
with rate $\kappa$ it makes a random walk jump and with rate $\gamma$ it flips velocity $v\to -v$.
If we denote $\mu(x,t,v)$ the probability to be at location $x\in\Z$ with velocity $v\in \{-1,1\}$ at time $t>0$, the generator
\eqref{gena} corresponds to the master equation (or Kolmogorov forward equation)
\beq\label{mastereqi}
\frac{d\mu(x,t,v)}{dt} &=& \lambda \mu(x-v,t,v) + \kappa (\mu(x-1,t, v)+ \mu(x+1,t,v)) + \gamma \mu(x,t,-v)
\nonumber\\
&-& (2\kappa+\lambda + \gamma) \mu(x,t,v)
\eeq
\subsection{Exact computation of the Fourier-Laplace transform}
The master equation \eqref{mastereqi} can be solved using Fourier-Laplace transform.
We define
\be\label{fourier}
\hat{\mu} (q,t,v)= \sum_{x} e^{iqx} \mu(x,t,v)
\ee
and view this quantity as a two-column, denoted $\overline{\mu}(q,t, \cdot)$ indexed by row index $v=1,-1$. The master equation \eqref{mastereqi} then becomes, after Fourier transform:
\be\label{bonk}
\frac{d}{dt} \overline{\mu}(q,t) = M(q) \overline{\mu}(q,t)
\ee
with $M(q)$ a symmetric two by two matrix of the form
\beq\label{M}
&&M(q)=
\left(
\begin{array}{cc}
a & b\\
b & a^*
\end{array}
\right)
\eeq
where $*$ denotes complex conjugate and where
\beq\label{abb}
a &=& (2\kappa +\lambda)(\cos(q)-1)-\gamma + i\lambda\sin(q) \nonumber\\
b &=& \gamma
\eeq
For the analysis of the scaling behavior of the position of the particle, it is convenient to further Laplace transform $\omu (q,t)$
i.e., we define, for $z>0$ the column vector
\be\label{lapma}
\homu(q,z)= \int_0^\infty \omu(q,t) e^{-zt}\ dt
\ee
then, from \eqref{bonk} we find
\[
\homu(q,z) = ( zI- M(q) )^{-1} {\bar{\mu}}_0 (q)
\]
For the initial position and velocity we choose $X_0=0$, and $v=\pm 1$ with probability $1/2$.
Then we have, ${\bar{\mu}}_0 (q)= \frac12 (1,1)^T$ where $T$ denotes transposition.
We further define the Fourier Laplace transform of the distribution of the  particle position:
\be
S(q,z)=\int_0^\infty \E e^{iqX_t} e^{-zt} \ dt = \sum_v \homu(q,z,v)=(1,1) \homu(q,z)
\ee

Then we have, using \eqref{lapma}
\be\label{sq}
S(q,z)= (\homu(q,z,1)) + (\homu(q,z, -1)) =
\frac12 (1,1) (zI-M(q))^{-1} (1,1)^T
\ee

Using the explicit formulas \eqref{M}, \eqref{abb}, we obtain
\be\label{sq}
S(q,z)= \frac{2\gamma + z - (\lambda +2\kappa) (\cos(q)-1)}{(\gamma +z -(\lambda+2\kappa) (\cos(q)-1))^2-\gamma^2 +  \lambda^2 \sin^2(q)}
\ee
For a more general velocity distribution at time zero, i.e., $X_0=0$, and $v=1$, resp.\ $v=-1$, with probability $\alpha$, resp. $1-\alpha$, we find
\be\label{genst}
S(q,z)= \frac{i\lambda (2\alpha -1) \sin(q)+ 2\gamma + z - (\lambda +2\kappa) (\cos(q)-1)}{(\gamma +z -(\lambda+2\kappa) (\cos(q)-1))^2-\gamma^2 +  \lambda^2 \sin^2(q)}
\ee
\subsection{Diffusive scaling behavior}
We can now use the explicit formula \eqref{sq} to obtain
the limit distribution of
$\epsi X_{\epsi^{-2} t}$ as $\epsi\to 0$. This amounts to
understand the scaling behavior of $\epsi^2 S(\epsi q, \epsi^2 z)$.
In particular  $\epsi X_{\epsi^{-2} t}\to \caN(0,\si^2 t)$ as $\epsi\to 0$ (in distribution), where $\caN(0,\si^2 t)$
denotes a normal with mean zero and variance $\si^2 t$, corresponds to the limiting
scaling behavior
\[
\lim_{\epsi\to 0} \epsi^2 S(\epsi q, \epsi^2 z)=  \frac{1}{z + \frac{q^2}{2}\si^2}
\]
If we obtain this scaling behavior, we call $\si^2$ the (limiting) diffusion constant.
We compute from the exact formula \eqref{sq}
\be\label{boci}
\lim_{\epsi\to 0} \epsi^2 S(\epsi q, \epsi^2 z) = \frac{1}{z + \frac{q^2}{2}\si^2}
\ee
with the limiting diffusion constant
\be\label{difcon}
\si^2= 2\kappa +\lambda + \frac{\lambda^2}{\gamma}
\ee
\br\label{bomom}
\bi
\item[a)] Notice that the limiting cases $\lambda=0$ corresponds to random walk at rate $\kappa$ to the left or right
which has diffusion constant $2\kappa$, and $\gamma\to\infty$ corresponds to a limiting Markovian random walk
for the position moving with rate $\kappa + \tfrac{\lambda}{2}$ to the right or to the left which has
diffusion constant $2\kappa +\lambda$. This is the so-called slow-fast limit.
Therefore, the extra term $\frac{\lambda^2}{\gamma}$ in \eqref{difcon} is the correction (w.r.t.\ slow-fast limit) to the diffusion constant due to the activity of the particles.
\item[b)] Using \eqref{genst}, we find that the scaling limit for a general initial velocity distribution is identical to the limit found in \eqref{boci}, \eqref{difcon}.
\item[c)]
If the random walk part of the generator is of a more general form than the nearest neighbor $\kappa (f(x+1, v)+ f(x-1,v) -2f(x,v))$ in \eqref{gena}, namely of the form
\[
\kappa\sum_{z\in\Z} p(z) (f(x+z,v)-f(x,v))
\]
where $p(z)=p(-z)$ is a symmetric transition rate
on $\Z$ such that $\sum_z z^2 p(z)<\infty$, then the same analysis applies, and the matrix $M(q)$ in \eqref{M} has the same off-diagonal element $b$, and the upper-diagonal element  $a$ has to be replaced by
\[
a=  \kappa \sum_{z} (\cos(qz)-1) p(z) + \lambda (\cos(q)-1) + i\lambda \sin(q)
\]
The scaling limit \eqref{boci} leads then to the diffusion constant
\be\label{gendifcon}
\si^2= \kappa \sum_z z^2 p(z) + \lambda + \frac{\lambda^2}{\gamma}
\ee
This shows that the addition to the diffusion constant due to activity is always equal to $\lambda + \frac{\lambda^2}{\gamma}$, irrespective of the precise form of the
random walk part of the generator.
\ei
\er
\subsection{Continuum limit: the telegrapher's process}
Let us now rescale the process $(X_t, v_t)$ with generator \eqref{gena} as follows: we consider the generator $L_\epsi$ on the state space
$\epsi\Z \times \{-1,1\}\subset \R \times \{-1,1\}$ defined by
\beq\label{epscontgen}
L_\epsi f(x,v) &=&
\epsi^{-1}\lambda (f(x+\epsi v, v) - f(x,v))
\nonumber\\
&+& \epsi^{-2}\kappa ( f(x+ \epsi, v) + f(x-\epsi,v)- 2 f(x,v))
\nonumber\\
&+& \gamma (f(x,-v)-f(x,v))
\eeq
This corresponds to a continuum limit in space, a diffusive rescaling of time and a rescaling of the parameters $\lambda\to \epsi\lambda$,
$\gamma\to\epsi^2\gamma$.

Then a Taylor
expansion gives
that for smooth functions $f: \R\times \{-1,1\}\to\R$ vanishing at $x\to\infty$, and for $(x_\epsi, v_\epsi)\in  \epsi\Z \times \{-1,1\}$ converging
to $(x,v)\in \R\times \{-1,1\}$  we obtain the limiting generator
$\lim_{\epsi\to 0}L_\epsi f(x_\epsi,v_\epsi) = \loc f(x,v)$ where
\be\label{contgen}
\loc  = \lambda v \partial_x + \kappa \partial^2_x + \gamma (\theta - I)
\ee
Here $I$ denotes the identity operator and $\theta f(x,v) = f(x,-v)$ is the velocity flip operator.

By the convergence of the generators, we obtain the weak convergence of the corresponding processes (in path space), i.e.,
\[
\{(\epsi X_{\epsi^{-2} t}, v_{\epsi^{-2}t})^{\epsi\lambda, \kappa, \epsi^2\gamma}, t\geq 0\}\to \{(\caX_t, \caV_t), t\geq 0\}
\]
as
$\epsi\to 0$.
\subsubsection{Fourier Laplace transform and diffusion limit in the telegrapher's process}
The limiting process
described by the generator $\loc$ is a diffusion-jump process where the particle drifts in the direction of the velocity,
with additional Brownian noise with variance $2\kappa$, and where the velocity flips according to a Poisson process with rate $\gamma$.
Let us denote by $\{(\caX_t, \caV_t), t\geq 0\}$ this limiting process with generator $\loc$. We call this
process the telegraphers process, abbreviated $TP(\lambda,\kappa,\gamma)$. This process is well-studied, see e.g. \cite{bombo,dhar, dhar2}.

Denote by $\mu(x,v,t)$ the probability density to find the particle at time $t$ at position $x$ with velocity $v$,
in the telegrapher's process $\{(\caX_t, \caV_t), t\geq 0\}$, then
we have the Kolmogorov forward equation
\beq
\partial_t \mu(x,t,v)
&=& -\lambda v \partial_x \mu(x,t,v) +\kappa \partial_x^2 \mu(x,t,v)
\nonumber\\
&+ & \gamma (\mu(x,t,-v)-\mu(x,t,v))
\eeq
As in \cite{bombo}, consider the Fourier transform $\omu(q,t)= \int e^{iqx} \mu(x, t, \cdot)$ viewed as a column indexed by the velocity $v\in \{-1,1\}$,
and $\homu (q,z)$ the Laplace transform of $\omu(q,t)$ w.r.t.\ $t$-variable, then we have
the  analogue
of \eqref{lapma}
\be\label{lapmac}
\homu (q,z)=( zI- \caM(q) )^{-1} \omu(q,0)
\ee
with
\be\label{MC}
( zI- \caM(q) )=
\left(
\begin{array}{cc}
z- i\lambda q + \kappa q^2 +\gamma & -\gamma\\
-\gamma & z- i\lambda q + \kappa q^2 +\gamma \\
\end{array}
\right)
\ee
This leads to an explicit formula for the Fourier-Laplace transform of the position, given starting position $\caX_0=0$, with starting velocity uniformly distributed on $\{-1,1\}$.
\beq\label{sqcont}
\caS(q,z) &=&\int_0^\infty\E\left(e^{iq\caX_t} e^{-zt}\ dt\right)\nonumber\\
&=&
\frac{2\gamma+ z+ \kappa q^2 }{ (z+\kappa q^2 +\gamma)^2 +\lambda^2 q^2 -\gamma^2}
\nonumber\\
\eeq
Taking the diffusion limit in the telegrapher's  process $TP(\lambda,\kappa,\gamma)$, then results in
\[
\lim_{\epsi\to 0} \epsi^2 \caS(\epsi q,\epsi^2 z)= \frac{1}{z + \si^2\frac{q^2}{2}}
\]
with the diffusion constant equal to
\be\label{contidif}
\si^2= 2\kappa + \frac{\lambda^2}{\gamma}
\ee
as found earlier in \cite{bombo}.

Compared with \eqref{difcon} we see that in the discrete case there is an additional term $\lambda$ in the diffusion constant
coming from the Poissonian noise for the jumps in the direction of the velocity, which is absent in the continuum limit (where
the motion in the direction of the velocity is purely deterministic).

In the process $TP(\lambda,\kappa,\gamma)$ by letting $\gamma\to\infty$, we obtain the slow-fast limit which gives
\[
\lim_{\gamma\to\infty}\caS(q,z)= \frac{1}{\kappa q^2 +z}
\]
which corresponds to a Brownian  process $B_{\sqrt{2}\kappa t}$, with variance $2\kappa t$.
\subsubsection{Adding an external field}
In the telegrapher's process, we can add drift corresponding to an external field $E$ by modifying the Kolmogorov forward equations
as follows
\beq\label{driftcont}
\partial_t \mu(x,t,v)
&=& -\lambda v \partial_x \mu(x,t,v) -2\kappa E \partial_x \mu(x,t,v)+ \kappa \partial_x^2 \mu(x,t,v)
\nonumber\\
&+& \gamma (\mu(x,t,-v)-\mu(x,t,v))
\eeq
We abbreviate this process $TP_E(\lambda,\kappa,\gamma)$, where $E$ refers to the external field.
\section{Large deviations}
\subsection{Lattice model}\label{freesec}
Instead of computing $S(q,z)$ in \eqref{sq} we can also directly compute the characteristic function
\be
C(q,t)= \E( e^{iqX_t} )
\ee
for the active particle on the lattice model with generator \eqref{gena}.
We choose the starting point $X_0=0$ and with random initial velocity, i.e., $v=\pm 1$ with probability $1/2$.
This amounts to compute the exponential of the matrix $M(q)$ from \eqref{M} which can be done using
diagonalization, and results in
\be
e^{t M(q)}= \frac{e^{tA}}{2\gamma B} G(t,q)
\ee
where $G(t,q)$ is given by the symmetric two by two matrix
\beq\label{tarzan}
&&G(t,q)=
\left(
\begin{array}{cc}
A_{11} & A_{12}\\
A_{12} & A_{11}^*
\end{array}
\right)
\nonumber\\
\eeq
where
\beq\label{bartaz}
A_{11} &=& -2\gamma\lambda i \sin(k) \sinh (Bt) + 2\gamma B \cosh (tB)
\nonumber\\
A_{12} &=& 2\gamma^2 \sinh(tB)
\eeq
and where
\beq
A &=& (\cos(k)-1)(2\kappa + \lambda) -\gamma
\nonumber\\
B &=& \sqrt{\gamma^2-\lambda^2 \sin^2(k)}
\eeq
This allows us to compute the moment generating function
via
\beq\label{momgen}
\E( e^{\alpha X_t}) & = & \frac12(1,1)e^{t M(-i \alpha)} (1,1)^T
\eeq
We can use the explicit formula to obtain the following large deviation result.
\bt
The position of the particle satisfies the large deviation principle, i.e.,
\be\label{boki}
\pee\left(\frac{X_t}{t} \approx x\right) \approx e^{-t I(x)}
\ee
The rate function $I$ is the Legendre transform of the large deviation free energy function $F$, i.e.,
\be\label{ratefunction}
I(x)= \sup_{\alpha\in \R} (\alpha x-F(\alpha))
\ee
where
\beq\label{free}
F(\alpha)
= (2\kappa +\lambda)(\cosh (\alpha) -1)
+   \sqrt{\gamma^2 + \lambda^2\sinh^2(\alpha)}-\gamma
\eeq
and where \eqref{boki} is shorthand for the large deviation principle, i.e.,
\beq
\limsup\frac1t\log \pee\left(\frac{X_t}{t}\in  B\right)\leq -\inf_{x\in B} I(x), \ B
\subset \R\ \text{closed}\nonumber\\
\liminf\frac1t\log \pee\left(\frac{X_t}{t}\in  G\right)\geq -\inf_{x\in G} I(x), \ G\subset \R\ \text{open}
\eeq
\et

\bpr
If we are interested in the large deviation properties of $X_t/t$ we compute the limiting cumulant generating function, or
large deviation free energy function, using \eqref{momgen}, \eqref{bartaz}, \eqref{tarzan}:
\beq\label{freei}
F(\alpha) = \lim_{t\to\infty} \frac1t\log\E \left(e^{\alpha X_t}\right)\nonumber\\
\eeq
Notice that, from \eqref{momgen}, it follows that $F(\alpha)$ in \eqref{free} is equal to the largest eigenvalue of the symmetric
matrix $M(-i\alpha)$ which is explicitly given by
\be\label{mia}
M(-i\alpha)=
\left(
\begin{array}{cc}
(2\kappa+\lambda) (\cosh(\alpha)-1) + \lambda\sinh(\alpha) -\gamma & \gamma\\
\gamma & (2\kappa+\lambda) (\cosh(\alpha)-1) - \lambda\sinh(\alpha)-\gamma
\end{array}
\right)
\ee
This gives
\beq
F(\alpha)= (2\kappa +\lambda)(\cosh (\alpha) -1)
+   \sqrt{\gamma^2 + \lambda^2\sinh^2(\alpha)}-\gamma
\eeq
From the computation of $F(\alpha)$, using the Gaertner-Ellis theorem \cite{dembo}, \cite{hol} we obtain the claimed large deviation principle.
\epr

Let us look at three relevant limiting cases for the ``free energy function'' $F$ from \eqref{free}.
\ben
\item[a)] Expanding the free energy function $F$ around $\alpha\approx 0$ gives
\[
F(\alpha)= \frac12D\alpha^2 + O(\alpha^4)
\]
with $D= 2\kappa + \lambda + \tfrac{\lambda^2}{\gamma}$.
This is consistent with the diffusion constant found in \eqref{difcon}.
The function $F(\alpha)$ in \eqref{free} can be analytically extended in a neighborhood of
the origin in the complex plane, and as a consequence, we can reobtain the central limit theorem (which we found via the scaling behavior of the characteristic function) from the large deviation
free energy, see \cite{bryc}.
\item[b)] In the limit $\gamma\to\infty$ the free energy function
becomes
\[
F(\alpha)= (\cosh(\alpha)-1)(2\kappa + \lambda)
\]
which corresponds to the large deviations of a symmetric random walk jumping
with rates $\kappa +\lambda/2$ to the right or left. This is indeed the (slow-fast) scaling limit of the process as we saw before. For large values of $\gamma$ we have
\[
F(\alpha)= (\cosh(\alpha)-1)(2\kappa + \lambda) + \frac{\lambda^2}{2\gamma} \sinh^2(\alpha)  + o(1/\gamma)
\]
Remark also that $F$ in \eqref{free} is non-increasing as a function of $\gamma$.
\item[c)] In the continuum limit we rescale $\lambda\to \epsi\lambda$, $\gamma\to \epsi^2\gamma$, $X_t\to \epsi X_{\epsi^{-2} t}$, we find
\be\label{dublim}
\lim_{\epsi\to 0}\lim_{t\to\infty}\frac1{t}\log  \E^{\epsi\lambda, \epsi^2\gamma}\left( e^{\alpha \epsi X_{\epsi^{-2} t}}\right) = \kappa\alpha^2 + \sqrt{\gamma^2+ \lambda^2 \alpha^2}-\gamma^2
\ee
which corresponds to the large deviation free energy of the continuum model (cf.\ subsection \ref{banka}, and see also \cite{bombo}), i.e., the limits $\epsi\to 0$ and
$t\to \infty$ in \eqref{dublim} commute.
\een
Finally, we point to a generalization as in remark \ref{bomom} above.
If the random walk part of the generator is of a more general form than the nearest neighbor $\kappa (f(x+1, v)+ f(x-1,v) -2f(x,v))$, namely of the form
\[
\kappa\sum_{z\in\Z} p(z) (f(x+z,v)-f(x,v))
\]
where $p(z)=p(-z)$ is a symmetric transition rate
on $\Z$ such that
\[
\Lambda(\alpha):= \sum_z (e^{\alpha z} -1)p(z)= \sum_{z} p(z) (\cosh(\alpha z)-1)<\infty
\]
then the same analysis applies and leads to the
large deviation free energy
\be\label{mgenf}
F(\alpha)= \lim_{t\to\infty} \log\E\left(e^{\alpha X_t}\right)= \Lambda(\alpha) + \lambda (\cosh(\alpha)-1) + \sqrt{\gamma^2 + \lambda^2\sinh^2(\alpha)}-\gamma
\ee
Indeed, the random walk part only changes the diagonal elements of the matrix in \eqref{mia} which now becomes
\[
\left(
\begin{array}{cc}
\Lambda(\alpha) +\lambda (\cosh(\alpha)-1) + \lambda\sinh(\alpha) -\gamma & \gamma\\
\gamma & \Lambda(\alpha) + \lambda (\cosh(\alpha)-1) - \lambda\sinh(\alpha)-\gamma
\end{array}
\right)
\]
where the term $(\cosh(\alpha)-1)(2\kappa)$ has been replaced by $\Lambda(\alpha)$. The largest eigenvalue of this matrix is given by \eqref{mgenf}. This implies that
the diffusion constant is given by
\[
\si^2= \lim_{t\to\infty} \frac1t\E(X^2_t) = F''(0)= \Lambda''(0) + \lambda + \frac{\lambda^2}{\gamma}
\]
which corresponds to the diffusion constant found earlier in \eqref{gendifcon}.
\subsection{Continuum model with drift and the Einstein relation}\label{banka}
We can also compute the free energy function corresponding to the large deviations
of the process in the telegrapher's process with drift $TP_E(\lambda,\kappa,\gamma)$, with master equation \eqref{driftcont}, via a similar diagonalization procedure.
This leads to
\be\label{FE}
F_E(\alpha)= \lim_{t\to\infty}\log \E e^{\alpha X_t} = \alpha^2\kappa + 2\alpha \kappa E + \sqrt{\lambda^2\alpha^2 + \gamma^2} - \gamma
\ee
as found earlier in \cite{bombo}.
We can then compute the asymptotic velocity of the particle:
\[
\lim_{t\to\infty} \frac1t \E(X_t) = F_E'(0)= 2\kappa E
\]
The limiting diffusion constant does not depend on $E$ and equals
\[
\lim_{t\to\infty} \frac1t \text{Var}(X_t) = F_E''(0)= 2\kappa + \frac{\lambda^2}{\gamma}
\]
As a consequence, the Einstein relation, relating the limiting velocity and the diffusion constant is violated as soon as $\lambda\not= 0$ and the the correction due to the activity is of
order $\lambda^2$.
\section{General models on $\Zd$ with finitely many velocities}
The fact that the large deviation free energy function $F$ can be computed as the largest eigenvalue of a symmetric matrix is true in much greater generality. The big advantage of the simple one-dimensional context is the simplicity of the explicit formulas.
In this section we sketch how to generalize the results.
The generalized (lattice) active particle model is a process $\{(X_t, v_t): t\geq 0\}$ with position $x\in\Zd$, and velocity $v$ taking values in a finite set: $v\in V\subset \Zd$. The generator reads
\be\label{gengen}
L = \lambda L_{t} + \kappa L_{d} +  \gamma L_f
\ee
where the three parts of the generator correspond to transport (i.e., motion in the direction of the velocity), diffusion (random motion) and flipping of the velocity, and are given by
\beq
L_t f(x,v) &=& f(x+v,v)-f(x,v)\nonumber\\
L_{d} f(x,v) &=& \sum_{z} p(z) (f(x+z,v)-f(x,v)) \nonumber\\
L_f f(x,v) &=& \sum_{v'\in V} \pi(v,v') (f(x,v')-f(x,v))
\eeq
Here $p(z)=p(-z)$ is a symmetric probability distribution on $\Zd$ such that
\beq\label{pasu}
\sum_{z\in\Zd} e^{\langle \alpha, z\rangle} p(z) = \sum_{z\in\Zd} \cosh\left(\langle \alpha, z\rangle\right)p(z) < \infty,\  \text{for all}\ \alpha\in\R^d
\eeq
 which represents the ``random walk'' jumps.
The assumption \eqref{pasu} is in order to be able to deal with large deviations for the particle position.
Furthermore, we assume that the velocity flip process transition rates $\pi(v,v')$ are such that they generate an irreducible continuous-time Markov chain on the finite set of velocities $V\subset\Zd$.

As a consequence of this assumption, we have a unique invariant measure for this velocity hop process. Let us  denote by $\{v_t, t\geq 0\}$ this velocity flip process which has generator
\be\label{velflipgen}
A f(v)=\sum_{v'\in V} \pi(v,v')( f(v')-f(v))
\ee
for functions $f: V\to\R$.
By the above stated assumptions, the process
$\{v_t,t\geq 0\}$ with generator $A$ satisfies occupation time large deviations, i.e., in the sense of the large deviation principle we have
\be\label{ldpoc}
\pee^{(A)} \left(\frac1T \int_0^T\delta_{v_s} \approx \mu\right) \approx e^{-T I_A(\mu)}
\ee
for $\mu$ a probability measure on $V$. Let us denote by $\caP(V)$ the set of probability measures on $V$. The rate function is given by the Donsker-Varadhan
formula
\be\label{dons}
I_A(\mu) = -\inf_{f>0} \int_V \frac{Af}{f} d\mu
\ee
\br
If $\pi(v,v')=\pi(v',v)$, then the uniform measure $\nu(v)= \frac{1}{|V|}$ on $V$ is the unique reversible measure, and the rate function
is given by the Dirichlet form \cite{dembo}
\be\label{iiaa}
I_A (\mu)= \left\langle\sqrt{ \frac{d\mu}{d\nu}}, - A \sqrt{\frac{d\mu}{d\nu}} \right\rangle_{L^2(\nu)} = -\gamma \sum_{v, v'}\sqrt{\mu(v)}\sqrt{\mu(v')} A_{v,v'}
\ee
\er

Let us further denote
\[
\Gamma (\alpha)= \sum_z p(z) (e^{\langle \alpha,  z\rangle}-1)=\sum_z p(z) (\cosh(\langle \alpha,  z\rangle)-1)
\]
where $\langle, \rangle$ denotes the Euclidean inner product. Further denote
\[
\caF(t,\alpha, v) = \sum_{x} \mu_t (x,v) e^{\langle \alpha, x \rangle}
\]
where we think of this as a column vector function of $\alpha, t$, where the column is indexed by $v\in V$.
Then we derive from the master equation corresponding to the generator \eqref{gengen} the equation
\be\label{caf}
\frac{d \caF(t,\alpha, \cdot)}{dt}= M(\alpha) \caF(t,\alpha, \cdot)
\ee
where $M(\alpha)$ is the symmetric matrix with diagonal entries given by
\be\label{pisa}
M(\alpha)_{vv}= \kappa\Gamma(\alpha) + \lambda( e^{\langle \alpha, v\rangle}-1) -\gamma=: \psi_\alpha(v) -\gamma
\ee
and off-diagonal entries
\[
M(\alpha)_{vv'}= \gamma\pi(v',v)= \gamma \pi(v,v') = M(\alpha)_{v'v}
\]
We then have the following theorem identifying the large deviation free energy function
as the largest eigenvalue of the matrix $M(\alpha)$. Moreover, by occupation time large deviations,  this eigenvalue in turn can be expressed
in a variational form, using Varadhan's lemma. This is the content of the following theorem.
\bt\label{gengenthm}
Let $\lambda(v,\alpha)$ denote the eigenvalues of the matrix $M(\alpha)$ (possibly degenerate). Let
\be\label{freegengen}
F(\alpha)= \lim_{t\to\infty} \frac1t\log\E( e^{\langle \alpha, X_t \rangle})
\ee
denote the large deviation free energy function for the position of the particle. Then we have
\beq\label{boranko}
F(\alpha) &=& \sup_{v\in V} \lambda(v, \alpha)
\nonumber\\
&=&\sup_{\mu\in \caP(V)} \left(\sum_v\psi_\alpha (v)\mu(v) -\gamma I_A(\mu)  \right)
\nonumber\\
&=&
\kappa\Gamma(\alpha) + \sup_{\mu\in \caP(V)} \left(\sum_v \lambda(e^{\langle\alpha, v\rangle}-1)\mu(v) -\gamma I_A(\mu)\right)
\eeq
where $\psi_\alpha$ is given in \eqref{pisa} and
where the rate function $I_A$ is given by \eqref{dons}.
In the symmetric case $\pi(v,v')=\pi(v',v)$ this specializes to
\[
F(\alpha)=\sup_{\mu\in \caP(V)} \left(\sum_v\psi_\alpha (v)\mu(v) +\gamma \sum_{v, v'} \sqrt{\mu(v)}\sqrt{\mu(v')} A_{v,v'}  \right)
\]
As a consequence, $X_t/t$ satisfies the large deviation principle with rate function $I(x)= \sup_\alpha (\langle\alpha, x\rangle - F(\alpha))$.
\et
\bpr
We give the proof in the symmetric case $\pi(v,v')=\pi(v',v)$, the general case is analogous (replacing the unitary diagonalization by a more general diagonalization).
From \eqref{caf} we have
\be\label{momgengen}
\E( e^{\langle \alpha, X_t \rangle})= \sum_v \caF(t,\alpha, v)= \langle \overline{1}, e^{tM(\alpha)} \mu_0 \rangle
\ee
where $\mu_0$ denotes column giving the distribution of initial velocities (we assume that the particle starts at the origin at time zero) and
$\overline{1}$ de column vector with all entries equal to 1. Let us denote by $e(v,\alpha), v\in V$ the orthonormal basis of eigenvectors of
the matrix $M(\alpha)$, with corresponding eigenvalues $\lambda(v,\alpha)$.
Then, choosing the initial velocity distribution $\mu_0$ uniform, we have
\[
\langle \overline{1}, e^{tM(\alpha)} \mu_0 \rangle= \sum_{v\in V} |\langle v,e(v,\alpha) \rangle|^2 e^{t\lambda(v, \alpha)}
\]
from which we conclude that
the large deviation free energy function
\be\label{freegengen}
F(\alpha)= \lim_{t\to\infty} \frac1t\log\E( e^{\langle \alpha, X_t \rangle})=\sup_{v\in V} \lambda(v, \alpha)
\ee
equals the largest eigenvalue of the matrix $M(\alpha)$, as we found in the two state velocity case before in subsection \ref{freesec} by explicit computation.
Notice now that the matrix $M(\alpha)$ is the sum of the Markov generator $\gamma A$ (where $A$ is defined in \eqref{velflipgen})
and a function of $v$ indexed by $\alpha$ (or a diagonal matrix), i.e.,
\be\label{fkfk}
M(\alpha)_{v,v'} =\gamma A_{v,v'} + \psi_\alpha (v) \delta_{v,v'}
\ee
where $A_{v,v'}= \pi(v,v') (1-\delta_{v,v'}) - \delta_{v,v'}\left(\sum_{v'\in V}\pi(v,v')\right)$ is the generator matrix of the velocity flip process, and where
$\psi_\alpha(v)=  \kappa
\Gamma(\alpha) + \lambda( e^{\langle \alpha, v\rangle}-1)$. We can obtain an alternative variational formula for the largest eigenvalue of $A + \psi_\alpha$, via the Feynman Kac formula.
Indeed, using the Feynman Kac formula,  we have from \eqref{fkfk}
\[
\left(e^{T M(\alpha)} f\right)(v)= \E^{(A)}_v \left(e^{\int_0^T \psi_\alpha (v_{s\gamma})ds} f(v_{T\gamma})\right)
\]
where $\E^{(A)}_v$ denotes expectation in the velocity process $\{ v_t,t\geq 0\}$ with generator $A$, starting from $v$ (notice that the process with generator
$\gamma A$ is then simply the time re-scaled process $\{ v_{t\gamma} : t\geq 0\}$).
As a consequence, we obtain for the largest eigenvalue of the matrix $M(\alpha)$ in \eqref{fkfk} the alternative formula
\[
\sup_{v\in V} \lambda(v, \alpha) = \lim_{T\to\infty} \frac1T\log\E^{(A)}\left( e^{\int_0^T \psi_\alpha (v_{s\gamma})ds}\right)
=\gamma\lim_{S\to\infty}\frac1S\log\E^{(A)}\left( e^{\frac1\gamma\int_0^S \psi_\alpha (v_{s})ds}\right)
\]
Then we use Varadhan's lemma, combined  with the occupation time large deviations
\eqref{ldpoc} with rate function \eqref{iiaa}, and obtain
\beq
F(\alpha)= \sup_{v\in V} \lambda(v, \alpha) &=& \gamma\lim_{T\to\infty} \frac1{T}\log\E^{(A)} e^{\frac1\gamma\int_0^{T} \psi_\alpha (v_s)ds}\nonumber\\
&= &\gamma\sup_{\mu\in \caP(V)} \left(\sum_v\frac1\gamma\psi_\alpha (v)\mu(v) - I_A(\mu)\right)\nonumber\\
&=& \sup_{\mu\in \caP(V)} \left(\sum_v\psi_\alpha (v)\mu(v) - \gamma I_A(\mu)\right)\nonumber\\
&=&
\sup_{\mu\in \caP(V)} \left(\sum_v\psi_\alpha (v)\mu(v) +\gamma \sum_{v, v'} \sqrt{\mu(v)}\sqrt{\mu(v')} A_{v,v'}  \right)
\nonumber
\eeq
Notice the third equality in \eqref{boranko} simply follows from the fact that $\Gamma(\alpha)$ does not depend on $v$. Therefore, we obtain \eqref{boranko}.
The corresponding large deviation result then follows via the Gaertner-Ellis theorem.
\epr
\br
\bi
\item[a)]
In the limit $\gamma\to\infty$, the pre-factor $\gamma $ in front of $I_A(\mu)$ forces $\mu$ to be equal to the unique stationary measure
$\nu$ of the velocity flip process, for which $I_A(\nu)=0$.
As a consequence, the large deviation free energy $F(\alpha)$ simply becomes the moment generating function of the random walk on $\Zd$ which
jumps from $x$ to $x+v$ with rate $\lambda \nu(v) + p(v)$ for $v\in V $, and from $x$ to $x+z$ with rate $p(z)$, for $z\not\in V$.
This is exactly the slow-fast limit.
\item[b)] The large deviation free energy function $F$ in \eqref{boranko} is a non-increasing function of $\gamma$. Indeed, for $\gamma'\geq \gamma$ we have, for all $\mu\in \caP(V)$
\[
\left(\sum_v\psi_\alpha (v)\mu(v) - \gamma' I_A(\mu)\right) \leq \left(\sum_v\psi_\alpha (v)\mu(v) - \gamma I_A(\mu)\right)
\]
As a consequence, the large deviation rate function (which is the Legendre transform of $F$) is a non-decreasing function of $\gamma$.
Since the rate function converges to the rate function of the slow-fast limit random walk, it follows that for finite $\gamma$, the rate function
is always smaller or equal than its slow-fast limit.
\ei
\er
{\bf Acknowledgement:}
The authors thank Gioia Carinci for useful discussions.

\end{document}